\documentclass[12pt]{article}

\usepackage{graphicx}

\usepackage{amsfonts}
\usepackage{amssymb}
\usepackage{latexsym}
\usepackage{amsmath,texdraw}

\usepackage[cp1251]{inputenc}
\usepackage[english]{babel}
\usepackage[T2A]{fontenc}
\usepackage{pgfplots}

\usepackage{color}

 \newtheorem{thm}{Theorem}[section]
 
 \newtheorem{lem}{Lemma}[section]
 \newtheorem{cor}{Corollary}[section]
 \newtheorem{exm}{Example}[section]
 
 \newtheorem{rem}{Remark}[section]

\newcommand{\sign}{\mathop{\rm sign }\medskip}
\newcommand{\llll}{{\left<\hskip -3pt\left<\right.\right.}\hskip -3pt}
\newcommand{\rrrr}{\hskip -3pt{\left.\left.\right>\hskip -3pt\right>\hskip -1pt}}
\newcommand{\D}{\displaystyle}

\newcommand{\rbx}{\hfill{\rule{1ex}{1ex}}}
\begin{document}

\vspace*{10mm}

\begin{center}
{\Large\textbf{A numerical point of view at 
the Gurov-Reshetnyak inequality on the real line}}
\end{center}

\vspace{5mm}

\begin{center}

\textbf{Victor D.~Didenko$^1$, Anatolii~A.~Korenovskyi$^2$,
          Nor~Jaidi~Tuah$^1$}

\vspace{2mm}

%    Only \author and \address are required; other information is
%    optional.  Remove any unused author tags.

%    author one information
$1.$ Universiti Brunei Darussalam, Bandar Seri Begawan, BE1410  Brunei,\\
  \texttt{diviol@gmail.com, norjaidi.tuah@ubd.edu.bn}

$2.$ Odessa I.I.~Mech\-ni\-kov Na\-tional Uni\-ver\-sity,
              Dvo\-ryan\-skaya 2, 65026 Odes\-sa, Ukraine, \texttt{anakor@paco.net}

 \end{center}

\vspace{5mm}

\begin{abstract}
A "norm" of power function in the Gurov-Reshetnyak class on the real
line is computed. Moreover, a lower bound for the norm
of the operator of even extension from the semi-axis to the whole
real line in the Gurov-Reshetnyak class is obtained from numerical
experiments.
 \end{abstract}

  \vspace{6mm}


\vspace{2mm}


\section{Introduction}
\label{intro}%Don't forget to give each section
Let us consider functions $f:R\mapsto\mathbb R^+$ where $R$ is an
interval of $\mathbb R$. In what follows,  $R$ is the real line
$\mathbb R$ or the semi-axis $\mathbb R^+=[0,\infty)$. We also
assume that the function $f$ is locally summable on $R$, i.e. it is
summable on each bounded subinterval $I$ of $R$.

The mean value of the function $f$ on a bounded interval $I$ is
defined by
\begin{equation*}
f_I=\frac1{|I|}\int_If(x)\,dx,
\end{equation*}
and the mean oscillation of this function is
\begin{equation*}
\Omega(f;I)=\frac1{|I|}\int_I\left|f(x)-f_I\right|\,dx,
\end{equation*}
where $|\,\cdot \,|$ denotes the Lebesgue measure.

For a given $\varepsilon\in(0,2]$ the Gurov--Reshetnyak class
$\mathcal{GR}=\mathcal{GR}(\varepsilon)=\mathcal{GR}_R(\varepsilon)$
is defined as the set of all non-negative functions  $f$ which are
locally summable on $R$ and such that the Gurov--Reshetnyak
condition
$$
\Omega(f;I)\le\varepsilon\,f_I
$$
is satisfied on all bounded intervals $I\subset R$~(see~\cite{GR}).
Note that since any non-negative function  $f$ on any
interval $I$ satisfies the inequality $\Omega(f;I)\le2f_I$,  the
class $\mathcal{GR}_R(2)$ is trivial and it coincides with the class
of all functions locally summable on $R$. However, if
$\varepsilon\in(0,2)$ then $\mathcal{GR}_R(\varepsilon)$ is a
non-trivial class (see \cite[P. 112]{K07}, \cite{KLS}). If $I$ is a
subinterval of $R$, then the expression $\left<\,f\,\right>_{I}=
{\Omega(f;I)}/{f_I}$ is called the relative oscillation of the
function $f$ on the interval $I$. Further, the term $\llll f\rrrr_R=
\sup\limits_{I\subset R}\left<\,f\,\right>_{I}$ is called the "norm"
of function $f$ in the Gurov-Reshetnyak class $\mathcal{GR}_R$.

One of the main properties of functions from the Gurov--Reshetnyak
class consists in the possibility to improve their summability
exponents. This property lays the foundation for numerous
applications of this class of functions. More precisely, for any
$\varepsilon\in (0,2)$ there are $p^+_{R}= p^+_{R}(\varepsilon)>1$
and $p^-_{R}= p^-_{R}(\varepsilon)<0$, such that the condition
$f\in\mathcal{GR}_R(\varepsilon)$ implies the local summability of
the function $f^p$ for any $p\in(p^-_{R},p^+_{R})$ (see.
  \cite{B}, \cite{F90}, \cite{F94}, \cite{GR}, \cite{F85}, \cite{I}, \cite{KLS},
\cite{W}). For $R=\mathbb R^+$, the exact limiting value
$p^+_{\mathbb R_+}=p^+_{\mathbb R_+}(\varepsilon)>1$ of the positive
summability exponent $p$ is the root of the equation
\begin{equation*}
\frac{p^p}{\left(p-1\right)^{p-1}}= \frac2\varepsilon,
\end{equation*}
and $p^-_{\mathbb R_+}=1-p^+_{\mathbb R_+}<0$. The sharpness of the
values $p^-_{\mathbb R_+}$ and $p^+_{\mathbb R_+}$ can be verified
by the use of the power functions $g(x)=x^{1/(p-1)}$ and
$h(x)=x^{-1/p}$ $(x\in\mathbb R_+,\,p>1)$, respectively. Thus
\begin{equation}\label{eq11}
\varepsilon_{\mathbb R_+}(p)\equiv \llll g\rrrr_{\mathbb R_+}= \llll
h\rrrr_{{\mathbb R_+}}= 2\frac{(p-1)^{p-1}}{p^p},
\end{equation}
\cite[pp. 131, 144]{K07}, \cite{K15}, \cite{K90},  \cite{K03a},
\cite{K04}. These examples also show that for functions
$f\in\mathcal{GR}_{\mathbb R_+}(\varepsilon)$, the function $f^p$ is
not necessarily locally summable in the limiting cases $p=p^-_{\mathbb
R_+}(\varepsilon)<0$ or $p=p^+_{\mathbb R_+}(\varepsilon)>1$.

On the other hand, for $R=\mathbb R$ the sharp limiting summability
exponents $p^-_{\mathbb R}(\varepsilon)<0$ and $p^+_{\mathbb
R}(\varepsilon)>1$ of functions $f\in\mathcal{GR}_{\mathbb
R}(\varepsilon)$ are not known. It is clear that $p^-_{\mathbb
R}(\varepsilon)\le p^-_{\mathbb R_+}(\varepsilon)$, $ p^+_{\mathbb
R}(\varepsilon)\ge p^+_{\mathbb R_+}(\varepsilon)$. Similarly to
$R=\mathbb R_+$, it is only natural to assume that for $R=\mathbb R$
the power functions $f_\alpha(x)=|x|^{\alpha}$ $(x\in\mathbb R,\
\alpha>-1)$ with $\alpha=1/(p-1)$ and $\alpha=-1/p$ $(p>1)$ are also
extremal ones. However, the computation of the corresponding
Gurov-Reshetnyak "norms" $\varepsilon^-_{\mathbb R}(p)\equiv \llll
f_{1/(p-1)}\rrrr_{\mathbb R}$ and $\varepsilon^+_{\mathbb
R}(p)\equiv\llll f_{-1/p}\rrrr_{\mathbb R}$ in this case is not as
simple as for $R=\mathbb R_+$. Nevertheless, it is shown in
\cite{K15} that $\varepsilon^-_{\mathbb R}(p)>\varepsilon_{\mathbb
R_+}(p)$ and $\varepsilon^+_{\mathbb R}(p)>\varepsilon_{\mathbb
R_+}(p)$.

One of the main results of the present work is the computation of
the "norm" $\llll f_\alpha\rrrr_{\mathbb R}$ of the function
$f_\alpha$ in the Gurov-Reshetnyak class on the real line $\mathbb
R$ (see Theorem~\ref{theo21} below). In particular, this theorem
implies the equation $\varepsilon^-_{\mathbb
R}(p)=\varepsilon^+_{\mathbb R}(p)\equiv\varepsilon_{\mathbb R}(p)$
$(p>1)$ (cf. Corollary~\ref{cor23}).

The above problem can be reformulated as follows: If a monotone
function $f$ belongs to the class $\mathcal{GR}_{\mathbb
R_+}(\varepsilon)$ for an $\varepsilon\in(0,2)$, then its even
extension to $\mathbb R$, which is also denoted by $f$, belongs to
the Gurov-Reshetnyak class $\mathcal{GR}_{\mathbb R}(\varepsilon')$
with an $\varepsilon'\in[\varepsilon,2)$ (see Lemma \ref{lem21}).
Therefore, one can also ask a question about the
norms\begin{equation*} \left\|\,{\bf
T}\,\right\|_{\mathcal{GR}}^{(\varepsilon)}\equiv\frac1\varepsilon
\sup\left\{\llll f\rrrr_{\mathbb R}:\
 \llll f\rrrr_{\mathbb R_+} =\varepsilon
\ (0<\varepsilon<2)\right\}, \  \left\|\,{\bf
T}\,\right\|_{\mathcal{GR}}\equiv
\sup_{0<\varepsilon<2}\left\|\,{\bf
T}\,\right\|_{\mathcal{GR}}^{(\varepsilon)}
\end{equation*}
of the operator ${\bf T}$ of the even extension of monotone functions
$f\in\mathcal{GR}_{\mathbb R_+}(\varepsilon)$ to the real line
$\mathbb R$.

In the Table~\ref{tab:rezults} below we report the results of
numerical calculations of the values of $\varepsilon_{\mathbb R}(p)$
for various $p>1$. Comparing these results with known values of
$\varepsilon_{{\mathbb R}_+}(p)$, one obtains lower bounds for the
norms $\left\|\,{\bf T}\,\right\|_{\mathcal{GR}}^{(\varepsilon)}$
and $\left\|\,{\bf T}\,\right\|_{\mathcal{GR}}$.

An analogous problem, concerning the norm of the operator of the
even extension for the class $BMO$ of monotone functions with
bounded mean oscillation, has been considered in \cite{Kl}. Some
estimates for the norm of such an extension have been obtained in
\cite{S}. It is remarkable that the lower estimate presented in
Remark~\ref{rem31} for the norm of the operator of the even
extension $\left\|\,{\bf T}\,\right\|_{\mathcal{GR}}$ of the present
work coincides with the lower estimate $\left\|\,{\bf
T}\,\right\|_{BMO}$ obtained in \cite{DKT} for the corresponding
operator of the even extension of monotone functions $f\in BMO$ (see
Remark \ref{rem32}).

\section{Gurov-Reshetnyak inequality for power functions on the real line}

Let us recall that the mean value $f_I=\gamma$ of function $f$ on a
subinterval $I$ is uniquely defined by the condition\footnote{By
$E(P)$ we denote the set of all points
 $x\in E$, satisfying the condition
$P=P(x)$.}
\begin{equation*}
\int_{I\left(f\ge \gamma\right)}\left(f(x)-\gamma\right)\,dx=
\int_{I\left(f\le\gamma\right)}\left(\gamma-f(x)\right)\,dx.
\end{equation*}
It is easily seen that
\begin{equation*}
\Omega(f;I)= \frac2{|I|}\int_{I\left(f\ge
f_I\right)}\left(f(x)-f_I\right)\,dx= \frac2{|I|}\int_{I\left(f\le
f_I\right)}\left(f_I-f(x)\right)\,dx.
\end{equation*}
According to \cite{K15}, the Gurov-Reshetnyak "norm" of any monotone
function $f$ on $\mathbb R_+$ can be computed by the formula
\begin{equation*}
\llll f\rrrr_{{\mathbb R}_+}= \sup_{b>0}\left<\, f\,\right>_{(0,b)},
\end{equation*}
and this fact will be used in what follows.

Let us show that the even extension of monotone functions from a
non-trivial Gurov-Reshetnyak class  $\mathcal{GR}_{\mathbb R_+}$
belong to a non-trivial class $\mathcal{GR}_{\mathbb R}$.

\begin{lem}\label{lem21}
For any $\varepsilon\in(0,2)$ there is an $\varepsilon'\in(0,2)$
such that if a monotone function $f$ belongs to
$\mathcal{GR}_{\mathbb R_+}(\varepsilon)$, then the even extension
of $f$ from ${\mathbb R_+}$ to $\mathbb R$ belongs to
$\mathcal{GR}_{\mathbb R}(\varepsilon')$.
\end{lem}

\textbf{Proof.}
The proof of this lemma can be split into three steps.

{\bf Step 1.} Consider an even function $f\in\mathcal{GR}_{\mathbb
R_+}(\varepsilon)$, $0<\varepsilon<2$. Then for any interval
$I\subset\mathbb R_+$ the Gehring inequality holds\footnote{Recall
the Gehring inequality first originated in \cite{G}.}
\begin{equation}\label{eq22}
\left(\frac1{|I|}\int_If^q(x)\,dx\right)^{1/q}\le
B\cdot\frac1{|I|}\int_If(x)\,dx,
\end{equation}
where $q>1$ and $B>1$ depend on the parameter  $\varepsilon$ only
\cite[P. 131]{K07}, \cite{K90}.

{\bf Step 2.} Let us show that the function $f$ satisfies the
Gehring inequality on $\mathbb R$. It suffices to consider only the
intervals of the form $I=(-a,b)$, where $0<a<b$. Since $f$ is an
even function, one can apply the inequality \eqref{eq22} on the
interval $(0,b)$ and obtain
\begin{multline*}
\left(\frac1{|I|}\int_If^q(x)\,dx\right)^{1/q}\le
\left(\frac2{a+b}\int_{0}^bf^q(x)\,dx\right)^{1/q}\le
\\
\le \left(\frac{2b}{a+b}\right)^{1/q}B\frac1b\int_{0}^bf(x)\,dx\le
2^{1/q}\left(\frac{b}{a+b}\right)^{1/q-1}B\frac1{a+b}\int_{-a}^bf(x)\,dx\le
\\
\le2B\frac1{|I|}\int_If(x)\,dx.
\end{multline*}

{\bf Step 3.} Now it remains to use the fact that the Gehring
inequality implies the Gurov-Reshetnyak inequality with an
$\varepsilon'\in(0,2)$ (see \cite[P. 114]{K07}, \cite{K03}), which
completes the proof.
\rbx

\begin{rem}\label{rem21}
The above obtained value of
$\varepsilon'\equiv\varepsilon'(\varepsilon)$ is not exact since at
each step of the proof of Lemma~\ref{lem21}, the parameters used are
overestimated.
\end{rem}

\begin{rem}\label{rem22}
For $\varepsilon<1$ there is a simpler proof of Lemma~\ref{lem21}.
In this case one can employ the simple inequality~\cite{DKT}
\begin{equation*}
\Omega(f;(-\delta b,b))\le\frac2{1+\delta}
\Omega(f;(0,b))\quad(0\le\delta\le1,\ b>0).
\end{equation*}
Indeed, since
\begin{equation*}
f_{(-\delta b,b)}=\frac{1}{(1+\delta)b}\int_{-\delta b}^bf(x)\,dx\ge
\frac1{1+\delta}\frac1b\int_0^bf(x)\,dx=
\frac1{1+\delta}f_{(0,b)},
\end{equation*}
then
\begin{equation*}
\frac{\Omega(f;(-\delta b,b))}{f_{(-\delta b,b)}}\le
(1+\delta)\frac{\Omega(f;(-\delta b,b))}{f_{(0,b)}}\le
(1+\delta)\frac{\frac2{1+\delta}\Omega(f;(0,b))}{f_{(0,b)}}\le
2\llll f\rrrr_{{\mathbb R}_+}.
\end{equation*}
Other details of the proof are left to the reader.
\end{rem}

\begin{rem}\label{rem23}
It follows from the proof in Remark~\ref{rem22} that if
$\varepsilon\in(0,1)$ then $\left\|\,{\bf
T}\,\right\|_{\mathcal{GR}}^{(\varepsilon)}\le2$. However, since
$\varepsilon'$ in Lemma~\ref{lem21} satisfies the relation
$\varepsilon'\le2$ the inequality $\left\|\,{\bf
T}\,\right\|_{\mathcal{GR}}^{(\varepsilon)}\le2$ remains valid for
$\varepsilon\in[1,2)$. Therefore, one also has $\left\|\,{\bf
T}\,\right\|_{\mathcal{GR}}\le2$. On the other hand, a lower bound
for $\left\|\,{\bf T}\,\right\|_{\mathcal{GR}}$, obtained in
numerical experiments, is presented in Remark~\ref{rem31} below.
\end{rem}

Further, we will compute the "norm" of power function in the
Gurov-Reshetnyak class. Using a linear transformation one can check
that for the function $f_\alpha(x)=|x|^\alpha$ $(x\in\mathbb
R,\alpha>-1)$ the following relations
\begin{equation*}
\llll f_\alpha\rrrr_{{\mathbb R}}=
\sup_{0\le\eta\le1}\left<\,f_\alpha\,\right>_{(-\eta,1)},\quad
\llll f_\alpha\rrrr_{{\mathbb R}_+}=
\left<\,f_\alpha\,\right>_{(0,1)}=
\frac{2|\alpha|}{(\alpha+1)^{(\alpha+1)/\alpha}}
\end{equation*}
hold.

\begin{thm}\label{theo21}
%{\bf Лемма I.}
%{\it
If $\alpha>-1$, $\alpha\ne0$, then
\begin{equation*}
\llll f_\alpha\rrrr_{{\mathbb R}}=
\llll f_\alpha\rrrr_{{\mathbb R_+}}\cdot
\max_{0\le\eta\le1}\psi(\alpha,\eta),
\end{equation*}
where
\begin{align*}
&\psi(\alpha,\eta)= \\[1ex]
 &=\left\{\begin{array}{ll}
\D\frac{\left(1+\eta^{\alpha+1}\right)^{1/\alpha}}{(1+\eta)^{(\alpha+1)/\alpha}}+
\frac{(\alpha+1)^{(\alpha+1)/\alpha}}\alpha
\left[\frac1{1+\eta^{\alpha+1}}-\frac1{1+\eta}\right],& \text{ if }\; 0\le\eta\le\eta_1,\\[10pt]
\D 2\frac{\left(1+\eta^{\alpha+1}\right)^{1/\alpha}}{(1+\eta)^{(\alpha+1)/\alpha}},& \text{ if }\; \eta_1\le\eta\le1,\\
\end{array}\right.
\end{align*}
and $\eta_1=\eta_1(\alpha)\in(0,1)$ is the root of the equation
\begin{equation}\label{eq23}
\eta^\alpha=\frac1{1+\alpha(\eta+1)}.
\end{equation}
%}
\end{thm}

\textbf{Proof.}
For any fixed  $\eta\in[0,1]$ we set $I= I(\eta)=(-\eta,1)$. Then
\begin{equation*}
\left(f_\alpha\right)_I=\frac1{1+\eta}\int_{-\eta}^1|x|^\alpha\,dx=
\frac{1}{\alpha+1}\cdot\frac{1+\eta^{\alpha+1}}{1+\eta}.
\end{equation*}
Let $\eta_1\in(0,1)$ be the root of the equation
$\eta=\left(\left(f_\alpha\right)_I\right)^{1/\alpha}$, cf.
\eqref{eq23}. It is easily seen that this equation is solvable and
the root is unique.

{\bf (a).} If $\eta\le\eta_1$, then
\begin{multline*}
\Omega(f_\alpha;I)=\frac{2\cdot\sign\alpha}{1+\eta}\int\limits_{\left(\left(f_\alpha\right)_I\right)^{1/\alpha}}^1\left(x^\alpha-\left(f_\alpha\right)_I\right)\,dx=
\\
%=\frac{2\cdot\sign\alpha}{1+\eta}\left[\frac1{\alpha+1}\left(1-\left(\left(f_\alpha\right)_I\right)^{(\alpha+1)/\alpha}\right)-
%f_I\left(1-\left(\left(f_\alpha\right)_I\right)^{1/\alpha}\right)\right]=
%\\
=\frac{2\cdot\sign\alpha}{1+\eta}\left[\frac\alpha{\alpha+1}\left(\left(f_\alpha\right)_I\right)^{(\alpha+1)/\alpha}-\left(f_\alpha\right)_I+\frac1{\alpha+1}\right].
\end{multline*}
Set
\begin{multline*}
\varphi_0(\alpha,\eta)\equiv\left<\,f_\alpha\,\right>_I=
\frac{\Omega(f_\alpha;I)}{\left(f_\alpha\right)_I}=
%\frac{2\cdot\sign\alpha}{1+\eta}\left[\frac\alpha{\alpha+1}\left(\left(f_\alpha\right)_I\right)^{1/\alpha}-1+\frac1{\alpha+1}\frac1{\left(f_\alpha\right)_I}\right]=
%\\
%= \frac{2\cdot\sign\alpha}{1+\eta}
%\left[\frac\alpha{\alpha+1}\left(\frac1{1+\eta}\frac{1}{\alpha+1}\left(1+\eta^{\alpha+1}\right)\right)^{1/\alpha}-1+
%\frac1{\alpha+1}\frac{(\alpha+1)(1+\eta)}{1+\eta^{\alpha+1}}\right]=
\\
=\frac2{1+\eta}\frac{|\alpha|}{(\alpha+1)^{(\alpha+1)/\alpha}}\left(\frac{1+\eta^{\alpha+1}}{1+\eta}\right)^{1/\alpha}-
\frac{2\cdot\sign\alpha}{1+\eta}+\frac{2\cdot\sign\alpha}{1+\eta^{\alpha+1}}.
\end{multline*}
Taking into account that
\begin{equation*}
\varphi_0(\alpha,0)= \left<\,f_\alpha\,\right>_{(0,1)}=
\frac{2|\alpha|}{(\alpha+1)^{(\alpha+1)/\alpha}},
\end{equation*}
one obtains
\begin{multline*}
\psi_0(\alpha,\eta)\equiv
\frac{\left<\,f_\alpha\,\right>_{I}}{\left<\,f_\alpha\,\right>_{(0,1)}}=
%\frac{\varphi_0(\alpha,\eta)}{\varphi_0(\alpha,0)}=
\\
=
\frac{\left(1+\eta^{\alpha+1}\right)^{1/\alpha}}{(1+\eta)^{(\alpha+1)/\alpha}}+
\frac{(\alpha+1)^{(\alpha+1)/\alpha}}\alpha
\left[\frac1{1+\eta^{\alpha+1}}-\frac1{1+\eta}\right].
\end{multline*}

{\bf (b).} On the other hand, if $\eta\ge\eta_1$, then
\begin{multline*}
\Omega(f_\alpha;I)=\frac{2\cdot\sign\alpha}{1+\eta}\cdot2\int\limits_0^{\left(\left(f_\alpha\right)_I\right)^{1/\alpha}}\left(\left(f_\alpha\right)_I-x^\alpha\right)\,dx=
\\
%=\frac{4\cdot\sign\alpha}{1+\eta}\left[\left(\left(f_\alpha\right)_I\right)^{(\alpha+1)/\alpha}-\frac1{\alpha+1}\left(\left(f_\alpha\right)_I\right)^{(\alpha+1)/\alpha}\right]=
=\frac4{1+\eta}\frac{|\alpha|}{\alpha+1}\left(\left(f_\alpha\right)_I\right)^{(\alpha+1)/\alpha}.
\end{multline*}
Consider now the expression
\begin{equation*}
\varphi_1(\alpha,\eta)\equiv\left<\,f_\alpha\,\right>_I=
\frac{\Omega(f_\alpha;I)}{\left(f_\alpha\right)_I}=
%\\
%=
%\frac4{1+\eta}\frac{|\alpha|}{\alpha+1}\left(\left(f_\alpha\right)_I\right)^{1/\alpha}=
%\frac4{1+\eta}\frac{|\alpha|}{\alpha+1}\left(\frac1{1+\eta}\frac1{\alpha+1}\left(1+\eta^{\alpha+1}\right)\right)^{1/\alpha}=
%\\
%=\frac4{(1+\eta)^{(\alpha+1)/\alpha}}\frac{|\alpha|}{(\alpha+1)^{(\alpha+1)/\alpha}}\left(1+\eta^{\alpha+1}\right)^{1/\alpha}
%=
\frac{4|\alpha|}{(\alpha+1)^{(\alpha+1)/\alpha}}\frac{\left(1+\eta^{\alpha+1}\right)^{1/\alpha}}{(1+\eta)^{(\alpha+1)/\alpha}}.
\end{equation*}
Since
\begin{equation*}
\varphi_1(\alpha,1)=
\frac{2|\alpha|}{(\alpha+1)^{(\alpha+1)/\alpha}}=
\varphi_0(\alpha,0)= \left<\,f_\alpha\,\right>_{(0,1)},
\end{equation*}
then
\begin{equation*}
\psi_1(\alpha,\eta)\equiv
\frac{\left<\,f_\alpha\,\right>_{I}}{\left<\,f_\alpha\,\right>_{(0,1)}}=
\frac{\varphi_1(\alpha,\eta)}{\varphi_1(\alpha,1)}=
2\frac{\left(1+\eta^{\alpha+1}\right)^{1/\alpha}}{(1+\eta)^{(\alpha+1)/\alpha}}.
\end{equation*}

Let $\psi$ be the function defined by
\begin{equation*}
\psi(\alpha,\eta)=\left\{\begin{array}{ll}
\psi_0(\alpha,\eta),\quad0\le\eta\le\eta_1,\\[1ex]
\psi_1(\alpha,\eta),\quad\eta_1\le\eta\le1.\\
\end{array}\right.
\end{equation*}
Then
\begin{equation*}
\frac{\llll f_\alpha\rrrr_{{\mathbb
R}}}{\llll f_\alpha\rrrr_{{\mathbb R_+}}}=
\max_{0\le\eta\le1}\psi(\alpha,\eta),
\end{equation*}
and the proof is completed.
\rbx

\begin{cor}\label{cor21}
Let $\psi$ be the function defined in Theorem~\ref{theo21}. Then
\begin{equation*}
\psi\left(-\frac\alpha{\alpha+1},\eta^{\alpha+1}\right)=
\psi\left(\alpha,\eta\right)\quad
(\alpha>-1,\ 0\le\eta\le1).
\end{equation*}
\end{cor}

\textbf{Proof.}
Straightforward calculations.
\rbx

For $p>1$ we set $\alpha=1/(p-1)$. Then
 $$\alpha+1= p/(p-1),\quad
-\alpha/(\alpha+1)=-1/p,
 $$
and Corollary~\ref{cor21} can be rewritten in the following form.

\begin{cor}\label{cor22}
Let $p>1$ and let $\psi$ be the function defined in
Theorem~\ref{theo21}. Then
\begin{equation*}
\psi\left(-\frac1p,\eta^{p/(p-1)}\right)=\psi\left(\frac1{p-1},\eta\right)\quad
(0\le\eta\le1).
\end{equation*}
\end{cor}

Recall that  the "norms" on the Gurov-Reshetnyak class are denoted
by $\varepsilon^-_{\mathbb R}(p)\equiv \llll
f_{1/(p-1)}\rrrr_{\mathbb R}$ and $\varepsilon^+_{\mathbb
R}(p)\equiv \llll f_{-1/p}\rrrr_{{\mathbb R}}$, $p>1$. However, one
has
\begin{equation*}
\llll f_{1/(p-1)}\rrrr_{{\mathbb R_+}}=
\llll f_{-1/p}\rrrr_{\mathbb R_+}\equiv
\varepsilon_{{\mathbb R}_+}(p),
\end{equation*}
and Theorem~\ref{theo21} and Corollary~\ref{cor22} lead to the
following result.

\begin{cor}\label{cor23}
If $p>1$, then
\begin{equation*}
\llll f_{1/(p-1)}\rrrr_{\mathbb R}=
\llll f_{-1/p}\rrrr_{\mathbb R}\equiv
\varepsilon_{{\mathbb R}}(p)=
\varepsilon_{{\mathbb R}_+}(p)\cdot
\max_{0\le\eta\le1}\psi\left(\frac1{p-1},\eta\right).
\end{equation*}
\end{cor}

The next corollary allows us to improve Theorem~\ref{theo21} by
using a better description of the set where the function $\psi$ can
attain its maximum.

\begin{cor}\label{cor24}
Let $\alpha>-1$ and $\alpha\ne0$. Then
\begin{equation*}
\max_{0\le\eta\le1}\psi(\alpha,\eta)=\max_{0\le\eta\le\eta_1}\psi_0(\alpha,\eta),
\end{equation*}
where the function $\psi_0(\alpha,\eta)$ is defined in the proof of
Theorem~\ref{theo21}, and the number $\eta_1=\eta_1(\alpha)\in(0,1)$
is derived from the equation~\eqref{eq23}.
\end{cor}

\textbf{Proof.}
Since the function $\psi(\alpha,\eta)$ is continuous on the interval
$[0,1]$ in the variable $\eta$, it suffices to show that for any
fixed  $\alpha$ the function $\psi_1(\alpha,\eta)$, defined in the
proof of Theorem~\ref{theo21}, is decreasing on the interval
$\left[\eta_1,1\right]$. For this we compute the derivative
\begin{equation*}
\frac\partial{\partial\eta}\psi_1(\alpha,\eta)=
%\\
%= \frac2{(1+\eta)^{2(\alpha+1)/\alpha}}
%\left[\frac1\alpha\left(1+\eta^{\alpha+1}\right)^{1/\alpha-1}(\alpha+1)\eta^\alpha(1+\eta)^{(\alpha+1)/\alpha}-
%\left(1+\eta^{\alpha+1}\right)^{1/\alpha}\frac{\alpha+1}\alpha(1+\eta)^{1/\alpha}\right]=
%\\
%=\frac2{(1+\eta)^{2(\alpha+1)/\alpha}}\frac{\alpha+1}\alpha(1+\eta)^{1/\alpha}\left(1+\eta^{\alpha+1}\right)^{1/\alpha-1}
%\left[\eta^\alpha(1+\eta)-\left(1+\eta^{\alpha+1}\right)\right]=
%\\
%=
2\frac{\alpha+1}\alpha
\frac{\left(1+\eta^{\alpha+1}\right)^{1/\alpha-1}}{(1+\eta)^{2+1/\alpha}}\left[\eta^\alpha
-1\right],
\end{equation*}
and observe that if $0<\eta<1$, then the derivative is negative.
This completes the proof.
\rbx

Corollaries~\ref{cor22} and~\ref{cor24} imply the following result:

\begin{cor}\label{cor25}
If $p>1$, then
\begin{equation*}
\max_{0\le\eta\le\eta_1\left(1/(p-1)\right)}\psi_0\left(\frac1{p-1},\eta\right)=
\max_{0\le\eta\le\eta_1\left(-1/p\right)}\psi_0\left(-\frac1{p},\eta\right).
\end{equation*}
\end{cor}

Let us compute the derivative of the function $\psi_0(\alpha,\eta)$.
Thus
\begin{align*}
\frac\partial{\partial\eta}\psi_0(\alpha,\eta)&=
\frac{(\alpha+1)^{(\alpha+1)/\alpha}}{\alpha}
\frac{\alpha+1}{1+\eta^{\alpha+1}}\frac1{1+\eta}
\\
&\quad\times
\left\{\left[1-\left(\frac1{\alpha+1}\frac{1+\eta^{\alpha+1}}{1+\eta}\right)^{1/\alpha}\right]
\left[\frac1{\alpha+1}\frac{1+\eta^{\alpha+1}}{1+\eta}-\eta^\alpha\right]
\right.
\\
&\quad-
\left[\frac\alpha{\alpha+1}\left(\frac1{\alpha+1}\frac{1+\eta^{\alpha+1}}{1+\eta}\right)^{(\alpha+1)/\alpha}\!-\!
\frac1{\alpha+1}\frac{1+\eta^{\alpha+1}}{1+\eta}+\frac1{\alpha+1}\right]
\\
&\phantom{\qquad-
\frac\alpha{\alpha+1}\left(\frac1{\alpha+1}\frac{1+\eta^{\alpha+1}}{1+\eta}\right)^{(\alpha+1)/\alpha}}
\left.
\times\;\eta^\alpha\frac{(\alpha+1)(1+\eta)}{1+\eta^{\alpha+1}}
\right\}.
 \end{align*}
Note that for $\eta=\eta_1$ the second factor in the second line of
the formula is equal to zero. Moreover, the expression in the third
line is nothing else but
$(1+\eta)\Omega(f_\alpha;I)/(2\sign\alpha)$, and
elementary computations lead to the following representation for the
derivative of the function $\psi_0(\alpha,\eta)$:
 \begin{multline*}
\frac\partial{\partial\eta}\psi_0(\alpha,\eta)=\frac{\alpha+1}\alpha\frac1{(1+\eta)^{2+1/\alpha}\left(1+\eta^{\alpha+1}\right)^2}
\times
\\
\times \left[
\left(1+\eta^{\alpha+1}\right)^{(\alpha+1)/\alpha}\left(\eta^\alpha-1\right)+
(\alpha+1)^{1/\alpha}\left(1+\eta^{\alpha+1}\right)^2(1+\eta)^{1/\alpha}-
\right.
\\
\left.
-
(\alpha+1)^{(\alpha+1)/\alpha}\eta^\alpha(1+\eta)^{2+1/\alpha}\right].
\end{multline*}
It is easily seen that
\begin{align*}
&\frac\partial{\partial\eta}\psi_0(\alpha,0)=
\frac{\alpha+1}\alpha\left[(\alpha+1)^{1/\alpha}-1\right]>0\quad
\text{if
}\alpha>0,\\[1ex]
&\frac\partial{\partial\eta}\psi_0(\alpha,0+)=+\infty\quad \text{if
} -1<\alpha<0.
\end{align*}
On the other hand,
\begin{equation*}
\frac\partial{\partial\eta}\psi_0\left(\alpha,\eta_1\right)=
-\frac{(\alpha+1)^{(3\alpha+1)/\alpha}}{2|\alpha|}
\frac{1+\eta_1}{\left(1+\eta_1^{\alpha+1}\right)^2}\eta_1^\alpha
\Omega\left(f_\alpha;I\left(\eta_1\right)\right)<0.
\end{equation*}
This leads to the following result.

\begin{cor}\label{cor26}
For each fixed $\alpha$ the function $\psi_0(\alpha,\eta)$  attains
its maximal value on the interval $\left[0,\eta_1\right]$ at the
inner point $\eta_{max}=\eta_{max}(\alpha)\in\left(0,\eta_1\right)$,
i.e. where
\begin{equation*}
\frac\partial{\partial\eta}\psi_0(\alpha,\eta)=0.
\end{equation*}
\end{cor}

Corollary~\ref{cor26} means that $\eta_{max}$ is the root of the
equation
\begin{multline}\label{eq35}
\left(1+\eta^{\alpha+1}\right)^{(\alpha+1)/\alpha}\left(\eta^\alpha-1\right)+
(\alpha+1)^{1/\alpha}\left(1+\eta^{\alpha+1}\right)^2(1+\eta)^{1/\alpha}-
\\
-
(\alpha+1)^{(\alpha+1)/\alpha}\eta^\alpha(1+\eta)^{2+1/\alpha}=0.
\end{multline}
However, even in the simplest case $\alpha=1$ the authors do not
know an analytic solution of this equation (see Example~\ref{ex31}
below).

\begin{rem}\label{rem24}
The numerical study of the behaviour of the function
$\psi_0(\alpha,\eta)$ for different values of $\alpha$ shows that
the derivative $\frac\partial{\partial\eta}\psi_0(\alpha,\eta)$ has
a unique root $\eta_{max}=\eta_{max}(\alpha)$ in the interval
$(0,\eta_1)$. Nevertheless, the authors do not know any rigorous
proof of this fact.
\end{rem}

\section{Numerical experiments, examples, and
comments}

Fix an $\varepsilon\in(0,2)$. Set
$p=p(\varepsilon)=p^+_{R}(\varepsilon)>1$,
$\alpha=\alpha(\varepsilon)=1/(p(\varepsilon)-1)$ and define
$\eta_1=\eta_1(\varepsilon)$ by~\eqref{eq23}. According to
Theorem~\ref{theo21} and Corollary~\ref{cor24}, one has
\begin{equation}\label{eq34}
\left\|\,{\bf T}\,\right\|_{\mathcal{GR}}^{(\varepsilon)}
\ge
\max_{0\le\eta\le\eta_1}\psi_0(\alpha,\eta)
\equiv C_\varepsilon,
\quad
\left\|\,{\bf T}\,\right\|_{\mathcal{GR}}
\ge
\sup_{0<\varepsilon<2}C_\varepsilon
\equiv C.
\end{equation}

Table \ref{tab:rezults} shows some values of the parameters
mentioned obtained in numerical experiments, where the
columns $6$ and $8$ contain the maximum points of the function
$\psi_0\left(\alpha,\eta\right)$ for $\alpha=1/(p-1)$ and  $\alpha
=-1/p$, correspondingly.
\begin{table}[h]
\caption{The Gurov-Reshetnyak "norms" and extremal points.
Numerical results.\label{tab:rezults}}
\begin{center}
\begin{tabular}{|r|l|l|l|l@{\hspace{-1mm}}|l|l|l@{\hspace{-1mm}}|}
\hline%\hline
$p\quad$  &
  $\varepsilon=\varepsilon_{\mathbb R_+}$ &
   $\varepsilon_{\mathbb R}\quad$  &
    $C_\varepsilon=\D\frac{\varepsilon_{\mathbb R}}{\varepsilon_{\mathbb R_+}}$    &
      $\D\alpha=\frac1{p-1}$   &
      $\quad \eta_{max}^+\ $  &
       $\D \alpha=-\frac1p $  &
         $\quad \eta_{max}^-\ $\\
\hline
\scriptsize{$1\hskip 15pt$}&
\scriptsize{$2\hskip 15pt$}&
\scriptsize{$3\hskip 15pt$}&
\scriptsize{$4\hskip 25pt$}&
\scriptsize{$5\hskip 15pt$}&
\scriptsize{$6\hskip 15pt$}&
\scriptsize{$7\hskip 19pt$}&
\scriptsize{$8\hskip 15pt$}\\[-2pt]
\hline\hline
   1.15 & 1.2813 & 1.4647 & 1.143133 & 6.6667 & 0.5484 &  -0.8696 & 0.0100\\
   \hline
   1.20 & 1.1647 & 1.3542 & 1.162679 & 5.0000 & 0.4936 &  -0.8333 & 0.0145\\
   \hline
   1.33 & 0.9493 & 1.1346 & 1.195193 & 3.0303 & 0.4030 &  -0.7519 & 0.0257\\
   \hline
   1.50 & 0.7698 & 0.9378 & 1.218204 & 2.0000 & 0.3372 &  -0.6667 & 0.0383\\
   \hline
   1.67 & 0.6513 & 0.8018 & 1.231116 & 1.4993 & 0.2982 &  -0.5999 & 0.0486\\
   \hline
   {\bf 2.00} & {\bf 0.5000} & {\bf 0.6224} & {\bf 1.244737} & {\bf 1.0000} & {\bf 0.2531} &  {\bf -0.5000} & {\bf 0.0640}\\
   \hline
   {\bf 3.00} & {\bf 0.2963} & {\bf 0.3726} & {\bf 1.257683} & {\bf 0.5000} & {\bf 0.2001} &  -0.3333 & 0.0895\\
   \hline
   6.00 & 0.1340 & 0.1692 & 1.263337 & 0.2000 & 0.1638 &  -0.1667 & 0.1141\\
   \hline
  11.00 & 0.0701 & 0.0886 & 1.264397 & 0.1000 & 0.1508 &  -0.0909 & 0.1248\\
   \hline
  21.00 & 0.0359 & 0.0454 & 1.264692 & 0.0500 & 0.1442 &  -0.0476 & 0.1309\\
   \hline
 101.00 & 0.0073 & 0.0093 & 1.264793 & 0.0100 & 0.1388 &  -0.0099 & 0.1361\\
   \hline
1001.00 & 0.0007 & 0.0009 & 1.264797 & 0.0010 & 0.1376 &  -0.0010 & 0.1373\\
   \hline
9999.00 & 0.0001 & 0.0001 & {\bf 1.264797} & 0.0001 & 0.1375 &  -0.0001 & 0.1374\\
   \hline
\end{tabular}
\end{center}
\end{table}

\begin{figure}[h]
%\centerline{\includegraphics[width=15cm]{pic-1}}

%%%%%%%%%%%%%% Fig. 1
%
%
%
\begin{center}
%\hskip -10mm
\begin{tikzpicture}
\pgfplotsset{width=6cm}
\begin{semilogxaxis} [
%\begin{axis} [
  title = The Gurov - Reshetnyak's "norms",
    xlabel = {$p$},
    legend pos = north east,
    ymin = 0,
    grid = major,
] \legend{
    $\varepsilon_{R_+}(p)$,
    $\varepsilon_{R}(p)$
};

\addplot[mark = otimes,
    mark options = {
        scale = 1.0,
        fill = pink,
        draw = black
    }] coordinates {

   (1.15,1.2813)
   (1.20,1.1647)
   (1.33,0.9493)
   (1.50,0.7698)
   (1.67,0.6513)
   (2.00,0.5000)
   (3.00,0.2963)
   (6.00,0.1340)
  (11.00,0.0701)
  (21.00,0.0359)
% (101.00,0.0073)
%(1001.00,0.0007)
%(9999.00,0.0001)
};

\addplot[mark = diamond,
    mark options = {
        scale = 1.3,
        fill = pink,
        draw = black
    }] coordinates {

   (1.15,1.4647)
   (1.20,1.3542)
   (1.33,1.1346)
   (1.50,0.9378)
   (1.67,0.8018)
   (2.00,0.6224)
   (3.00,0.3726)
   (6.00,0.1692)
  (11.00,0.0886)
  (21.00,0.0454)
% (101.00,0.0093)
%(1001.00,0.0009)
%(9999.00,0.0001)

};

%\addplot{ 2*((x-1)^(x-1))/(x^x) };

\end{semilogxaxis}
\end{tikzpicture}
%
%%%%%%%%%%%%%%%%%%%% Fig 2
%
%
%
%
%\hskip 10mm
\begin{tikzpicture}
\pgfplotsset{width=6cm}
\begin{semilogyaxis} [
%\begin{axis} [
title = The limiting exponent,
    xlabel = {$\varepsilon$},
    legend pos = north east,
    ymin = 1.14,
    grid = major
] \legend{
   $p\left(\varepsilon_{R_+}\right)$,
   $p\left(\varepsilon_{R}\right)$
};

\addplot[mark = otimes,
    mark options = {
        scale = 1.0,
        fill = pink,
        draw = black
    }] coordinates {

   (1.2813,1.15)
   (1.1647,1.20)
   (0.9493,1.33)
   (0.7698,1.50)
   (0.6513,1.67)
   (0.5000,2.00)
   (0.2963,3.00)
   (0.1340,6.00)
  (0.0701,11.00)
  (0.0359,21.00)
% (0.0073,101.00)
%(0.0007,1001.00)
%(0.0001,9999.00)

 };

\addplot[mark = diamond,
    mark options = {
        scale = 1.3,
        fill = pink,
        draw = black
    }] coordinates {

   (1.4647,1.15)
   (1.3542,1.20)
   (1.1346,1.33)
   (0.9378,1.50)
   (0.8018,1.67)
   (0.6224,2.00)
   (0.3726,3.00)
   (0.1692,6.00)
  (0.0886,11.00)
  (0.0454,21.00)
% (0.0093,101.00)
%(0.0009,1001.00)
%(0.0001,9999.00)

};

\end{semilogyaxis}
\end{tikzpicture}
\end{center}
\centering
 \vskip -0.5cm \caption{Relations between $p$ and
$\varepsilon$.} \label{ris:pic-1}
\end{figure}
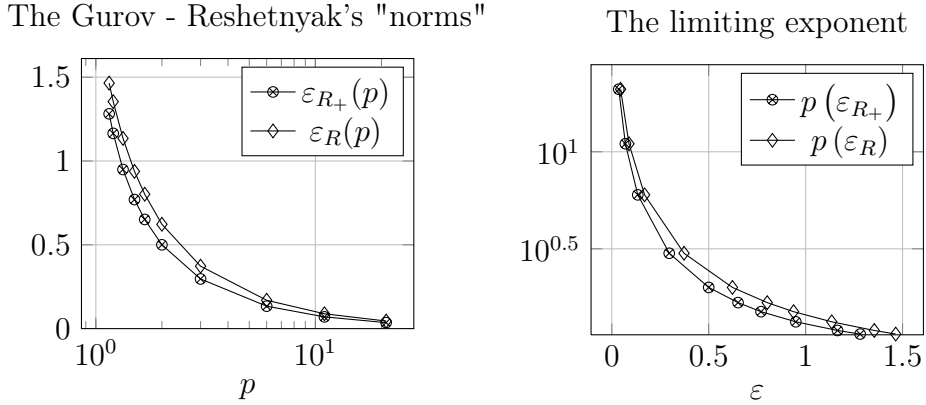

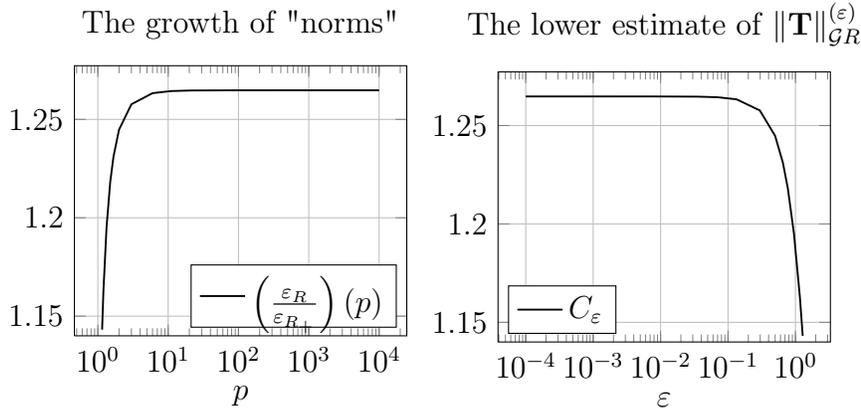
\begin{figure}[h]
%\centerline{\includegraphics[width=15cm]{pic-2}}

\begin{center}
%\hskip -10mm

%%%%%%%%%%%%%% Fig. 3
%
%\hskip -10mm
\begin{tikzpicture}
\pgfplotsset{width=6cm}
\begin{semilogxaxis} [
%\begin{axis} [
title = The growth of "norms",
    xlabel = {$p$},
    legend pos = south east,
    ymin = 1.14,
    grid = major
] \legend{
   $\left(\frac{\varepsilon_{R}}{\varepsilon_{R_+}}\right)(p)$
};

\addplot [solid, draw = black, line width=0.7pt]coordinates {

   (1.15,1.143133)
   (1.20,1.162679)
   (1.33,1.195193)
   (1.50,1.218204)
   (1.67,1.231116)
   (2.00,1.244737)
   (3.00,1.257683)
   (6.00,1.263337)
  (11.00,1.264397)
  (21.00,1.264692)
 (101.00,1.264793)
(1001.00,1.264797) (9999.00,1.264797)

};

\end{semilogxaxis}
\end{tikzpicture}
%
%
%
%%%%%%%%%%%%%%%%%%%%% Fig. 4
%
%
%
%\hskip 10mm
\begin{tikzpicture}
\pgfplotsset{width=6cm}
\begin{semilogxaxis} [
%\begin{axis} [
title = The lower estimate of $\|{\bf T}\|_{{\mathcal GR}}^{(\varepsilon)}$,
    xlabel = {$\varepsilon$},
    legend pos = south west,
    ymin = 1.14,
    grid = major
] \legend{
   $C_\varepsilon$
};~\addplot [solid, draw = black, line width=0.7pt]coordinates {

   (1.2813,1.143133)
   (1.1647,1.162679)
   (0.9493,1.195193)
   (0.7698,1.218204)
   (0.6513,1.231116)
   (0.5000,1.244737)
   (0.2963,1.257683)
   (0.1340,1.263337)
   (0.0701,1.264397)
   (0.0359,1.264692)
   (0.0073,1.264793)
   (0.0007,1.264797)
   (0.0001,1.264797)

};

\end{semilogxaxis}
\end{tikzpicture}
\end{center}

\vskip -0.5cm \caption{The growth of the norms during the extension
from $\mathbb R_+$ to $\mathbb R$.} \label{ris:pic-2}
%\centerline{Графики функций $\psi^+(1,\eta)$ и $\psi^-\left(-\frac12,\eta\right)$ $(p=2)$.}
\end{figure}

In addition, the results of numerical experiments are reflected in
Figures~\ref{ris:pic-1} and~\ref{ris:pic-2}.

%\begin{figure}[h]
%\centerline{\includegraphics[width=15cm]{pic-1}}
%\caption{Соотношения между $p$ и $\varepsilon$.}
%\label{ris:pic-1}
%\end{figure}

%\begin{figure}[h]
%\centerline{\includegraphics[width=15cm]{pic-2}}
%\caption{Рост ''норм'' при переходе от $\mathbb R_+$ к $\mathbb R$.}
%\label{ris:pic-2}
%\end{figure}

{\bf Comments to the graphs.}\nopagebreak

{\bf Figure \ref{ris:pic-1}.} The lower graph in the left part of
the Figure~\ref{ris:pic-1} shows the dependance of the
Gurov-Reshetnyak  "norms" $\varepsilon_{\mathbb R_+}(p)$ on the
parameter $p$. These results are obtained from formula~\eqref{eq11}.
Note that the data are represented in the logarithmic scale and do
not include all results from Column $2$. The upper graph shows the
dependance of the exponents $\varepsilon_{\mathbb R}(p)$ on $p$,
presented in Corollary~\ref{cor23}.

In the right part of Figure~\ref{ris:pic-1} we show the graphs of
the inverse relations, i.e. these are values of those $p$, for which
the function $f_{1/(p-1)}$ belongs to the class
$\mathcal{GR}_{\mathbb R_+}(\varepsilon)$ (the lower line) or to the
class  $\mathcal{GR}_{\mathbb R}(\varepsilon)$ (the upper line) for
a given $\varepsilon$.

{\bf Figure \ref{ris:pic-2}.} In the left graph we show the quotient
of the "norms" of the function $f_{1/(p-1)}$ in classes
$\mathcal{GR}_{\mathbb R}$ and $\mathcal{GR}_{\mathbb R_+}$ with
respect to the parameter~$p$. The right graph reflects the behaviour
of the parameter $C_\varepsilon$ (see~\eqref{eq34}) with respect to
$\varepsilon$. In both cases the variables $p$ and $\varepsilon$ are
shown in the logarithmic scale.

\begin{rem}\label{rem31}
As is seen from Table~\ref{tab:rezults}, the "norm" of the operator
$\mathbf{T}$ can be estimated as follows
 $$
\left\|\,{\bf T}\,\right\|_{\mathcal{GR}}\ge
C=\lim_{\varepsilon\to0+}C_\varepsilon\approx 1.264797,
 $$
where the constant $C$ is defined in~\eqref{eq34}. The corresponding
numerical value $C$ is shown in boldface in the last row of
Table~\ref{tab:rezults}.
\end{rem}

\begin{rem}\label{rem32}
For functions $f\in BMO$ with bounded mean oscillation the question
about the sharp value of the norm $\|\,f\,\|_{BMO,{\mathbb
R}}=\sup\limits_{I\subset {\mathbb R}}\Omega(f;I)$ of the even
extension from ${\mathbb R_+}$ to ${\mathbb R}$ of a
monotone function on the semi-axis $\mathbb R_+$ is posted in
\cite{Kl} and, to the best of our knowledge, it is still open. In
\cite{DKT} the $BMO$-norm of the function $f_0(x)=\ln(1/|x|)$, which
is a typical representative of this class, is found. Thus
\begin{equation*}
\left\|\,f_0\,\right\|_{BMO,\mathbb R}=
\frac2{\rm e}\cdot\frac1{t+1}\left[\exp\left(\frac{t\,\ln t}{t+1}\right)+ {\rm
e}\frac{t\,\ln t}{t+1}\right],
\end{equation*}
where $t>1$ is the root of the equation
\begin{equation}\label{eq36}
\exp\left(\frac{t\,\ln t}{t+1}\right)= {\rm
e}\left(t-1-\frac{t+1}{\ln t}\right).
\end{equation}
Since $\left\|\,f_0\,\right\|_{BMO,\mathbb R_+}=2/{\rm e}$, the
approximate solution of the equation \eqref{eq36} leads to the
following estimate (see \cite{DKT})
\begin{equation*}
\left\|\,{\bf T}\,\right\|_{BMO}=\!\!\sup_{\scriptsize
\begin{array}{cc}f\,\text{ is even on }\mathbb R,\\ \text{and monotone on }\mathbb
R_+\end{array}}\!\!\! \frac{\|\,f\,\|_{BMO,\mathbb
R}}{\|\,f\,\|_{BMO,\mathbb R_+}} \ge\frac{\|\,f_0\,\|_{BMO,\mathbb
R}}{\|\,f_0\,\|_{BMO,\mathbb R_+}}\equiv C_0 \approx 1.264797.
\end{equation*}

It is remarkable that $C$ and $C_0$ coincide up to $6$ significant
digits after the point, i.e. the lower bounds of the norms
$\left\|\,{\bf T}\,\right\|_{\mathcal{GR}}$ and $\left\|\,{\bf
T}\,\right\|_{BMO}$ of the operator of even extension ${\bf T}$
turned out to be the same within the calculation accuracy.
\end{rem}

\begin{exm}\label{ex31}
In the simplest case $\varepsilon=1/2$, $p=2$ and $\alpha=1$ or
$\alpha=-1/2$, one can find the roots of the equation~\eqref{eq23},
which are $\eta_1(1)=\sqrt2-1\approx0.414$,
$\eta_1\left(-1/2\right)=3-2\sqrt2\approx0.172$.
\begin{figure}[h]
\centering
\includegraphics%
[width=11cm]{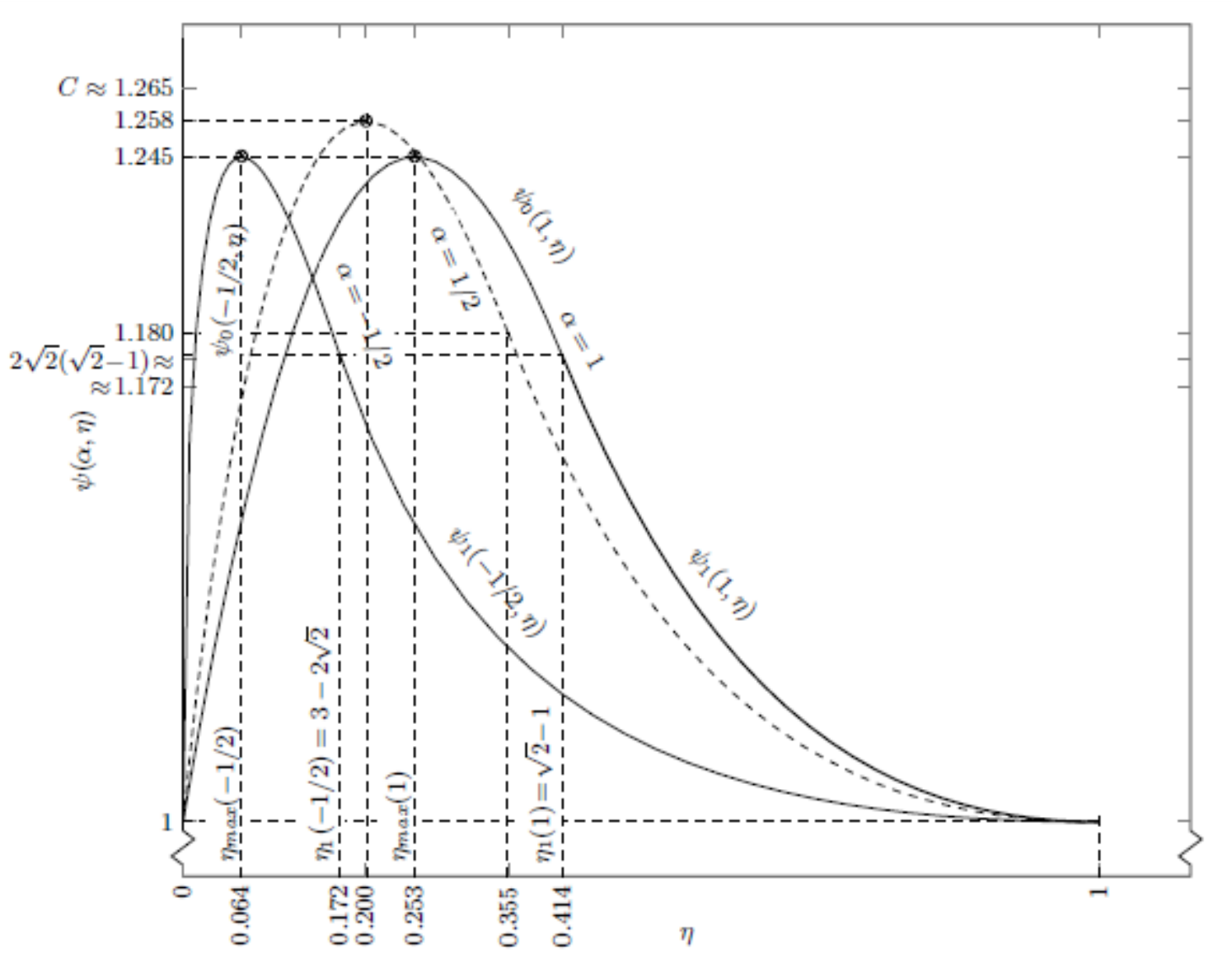} \caption{The graphs of the functions
$\psi(1,\eta)$, $\psi\left(-1/2,\eta\right)$ $(p=2)$ and
$\psi\left(1/2,\eta\right)$ $(p=3)$.} \label{ris:graph1}
\end{figure}

The corresponding line in Table~\ref{tab:rezults}  is marked out by
boldface. In Figure~\ref{ris:graph1} the graphs of the functions
$\psi\left(-1/2,\eta\right)$ and $\psi(1,\eta)$ are represented by
solid lines. Note that $\psi_0(1,\eta_1(1))=
\psi_0\left(-1/2,\eta_1\left(-1/2\right)\right)=
2\sqrt2\left(\sqrt2-1\right)\approx1.172$. For these values of the
parameters, the equation~\eqref{eq35}, which is used to find
$\eta_{max}=\eta_{max}(1)\in\left(0,\eta_1\right)$, takes the form
\begin{equation*}
3\eta^5-3\eta^4-6\eta^3-10\eta^2-\eta+1=0.
\end{equation*}
The detection of any analytic solution of this equation seems to be
difficult. However, one can show that in the interval
$\left(0,\eta_1\right)$ this equation has a unique solution  (see
Remark~\ref{rem24}). Indeed, since the second derivative of the
function $\psi_0(1,\eta)$ is
\begin{equation*}
\frac{\partial^2}{\partial\eta^2}\psi_0(1,\eta)=
-4\left(\frac{3\eta}{(1+\eta)^4}+2\frac{1-3\eta^2}{\left(1+\eta^2\right)^3}\right)
\end{equation*}
and $1-3\eta^2>0$ for $0<\eta<\sqrt2-1$, one obtains that
$\frac{\partial^2}{\partial\eta^2}\psi_0(1,\eta)<0$ on the interval
$\left(0,\eta_1\right)$. This means that
$\frac{\partial}{\partial\eta}\psi_0(1,\eta)$ is strictly decreasing
in the interval $\left(0,\eta_1\right)$, hence it has a unique root
in this interval.  By Corollary~\ref{cor22},  the derivative
$\frac{\partial}{\partial\eta}\psi_0(-1/2,\eta)$ also has a unique
root in the interval $\left(0,\eta_1(-1/2)\right)$.

Notice that the values $\eta_{max}(1)\approx0.253$,
$\eta_{max}\left(-1/2\right)\approx0.064$ are approximations of the
corresponding roots and $C_{1/2}\approx1.245$.
\end{exm}

\begin{exm}\label{ex32}
Let $\alpha=1/2$, $p=3$, and $\varepsilon\approx0.296$. In
Table~\ref{tab:rezults} the numerical results corresponding to this
case are shown in the line partially written in boldface. In
Figure~\ref{ris:graph1} the graph of the corresponding function
$\psi(1/2,\eta)$ is represented by the dashed line.

In this case, the equation~\eqref{eq23} takes the form
\begin{equation*}
\sqrt\eta=\frac2{3+\eta}.
\end{equation*}
The solution of this equation
\begin{equation*}
\eta_1\equiv\eta_1\left(\frac12\right)=\sqrt[3]{3+2\sqrt2}+\sqrt[3]{3-2\sqrt2}-2
\approx0.355
\end{equation*}
is obtained by Cardano formulas. We also have
\begin{equation*}
\psi_0\left(\frac12,\eta\right)=
\frac{\left(1+\eta^{3/2}\right)^2}{(1+\eta)^3}+
\frac{27}4\left[\frac1{1+\eta^{3/2}}-\frac1{1+\eta}\right],
\end{equation*}
$\psi\left(1/2,\eta_1\right)\approx1.180$. In addition, the equation
~\eqref{eq35} takes the form
\begin{equation*}
\left(1+\eta^{3/2}\right)^3\left(\eta^{1/2}-1\right)+
\frac94\left(1+\eta^{3/2}\right)^2(1+\eta)^2-
\frac{27}8\eta^{1/2}(1+\eta)^4=0,
\end{equation*}
and its approximate solution is $\eta_{max}(1/2)\approx0.200$.
Finally, we obtain an approximate value of $C_{0.296}$, namely
$C_{0.296}\approx1.258$.
\end{exm}

\section{Acknowledgements}
This research was supported by the Universiti Brunei Darussalam
under Grant UBD/GSR/S\&T/19.

\end{document}